\documentclass[12pt,reqno]{amsart}

\newcommand\version{March 26, 2010}

\usepackage{amsmath,amsfonts,amsthm,amssymb,amsxtra}

\newtheorem{theorem}{Theorem}[section]

\newtheorem{lemma}[theorem]{Lemma}

\theoremstyle{definition}

\theoremstyle{remark}

\numberwithin{equation}{section}

\newcommand{\C}{\mathbb{C}}

\renewcommand{\epsilon}{\varepsilon}

\newcommand{\N}{\mathbb{N}}

\renewcommand{\phi}{\varphi}
\newcommand{\R}{\mathbb{R}}

\newcommand{\Sph}{\mathbb{S}}

\newcommand{\lanbox}{\hfill \hbox{$\, \vrule height 0.25cm width 0.25cm depth 0.01cm \,$}}

\begin{document}

\title[Spherical positivity --- \version]{Spherical reflection positivity and \\ the Hardy-Littlewood-Sobolev inequality}

\author{Rupert L. Frank}
\address{Rupert L. Frank, Department of Mathematics,
Princeton University, Washington Road, Princeton, NJ 08544, USA}
\email{rlfrank@math.princeton.edu}

\author{Elliott H. Lieb}
\address{Elliott H. Lieb, Departments of Mathematics and Physics,
Princeton University,
       P.~O.~Box 708, Princeton, NJ 08544, USA}
\email{lieb@princeton.edu}

\thanks{\copyright\, 2009 by the authors. This paper may be reproduced, in its entirety, for non-commercial purposes.\\
Support through U.S. NSF grant PHY 0652854 is gratefully acknowledged.}

\begin{abstract}
We introduce the concept of spherical (as distinguished from planar) reflection positivity and use it to obtain a new proof of the sharp constants in certain cases of the HLS and the logarithmic HLS inequality. Our proofs relies on an extension of a work by Li and Zhu which characterizes the minimizing functions of the type $(1+|x|^2)^{-p}$.
\end{abstract}

\maketitle


\section{Introduction}

The well-known functions on $\R^N$, $f(x)=c(b^2+|x-a|^2)^{-p}$, where $a\in\R^N$, $b>0$ and $c>0$, appear as the optimizers in some classical functional inequalities, notably the Hardy-Littlewood-Sobolev (HLS) inequality and its dual, the Sobolev inequality. Given their ubiquity, these functions must be endowed with some special property, and this was identified by Y. Y. Li and M. J. Zhu \cite{LZh}. It is the property of reflection invariance through spheres as we shall explain later. One of our contributions is the proof, using reflection positivity through spheres, that optimizing functions for the HLS inequality must have this interesting reflection property. In this note we shall explain the classical reflection positivity through planes, the new reflection positivity through spheres and its application to the HLS inequality. This solves the problem of the sharp constant in this inequality \emph{without using symmetric decreasing rearrangement}, as was done earlier \cite{Li}. This note is a summary of our results in \cite{FrLi}, but contains partially alternate proofs of some topics, notably the proof of spherical reflection positivity, utilizing Gegenbauer polynomials instead of relying on conformal invariance. We also extend our analysis to the logarithmic version of the HLS inequality.

 Since we do not know, a priori, that our optimizing functions have the required continuity property needed for direct application of the Li--Zhu characterization lemma, we extend their lemma from functions to measures. The proof of this extension is given in \cite{FrLi}.

\subsection{Reflection positivity}

We begin with {\bf reflection positivity through} \textbf{planes} in $\R^N$. It is beloved of quantum field theorists \cite{OsSc,GlJa} and statistical mechanicians \cite{Li0,FrIsLiSi}, but it surely must have been
known to potential theorists in the nineteenth century. Consider the plane in $\R^N$ defined
by $x_N =0$ and a function $f $ with support in the half-space $H=\{x: x_N>0\} $ as well as its
reflected function $\theta f$ with support in the complementary half-space $\theta H =\{x: x_N<0\}$,
defined by $\theta f(x', x_N) = f(x', -x_N)$. Then RP states that
\begin{equation}     \label{eq:rp}
  I_\lambda[\theta f,f] := \iint_{\R^N\times\R^N} \frac{\overline{\theta f(y)}f(x)}{|y-x|^{\lambda}} \,dx\,dy \geq 0 \qquad {\rm for} \ N> \lambda\geq N-2\ {\rm and } \ \lambda >0.
\end{equation}
It is important to note that $f$ does not have to be positive, or even real.

The case $\lambda = N-2$ and $N\geq 3$ is, of course, the kernel of ordinary 
potential theory (the inverse of the Laplacian) -- and this is the classically known case of RP. The other cases may not have been known and are proved in our paper \cite{FrLi}, although an inequality equivalent to \eqref{eq:rp} was proved in \cite{LoMa}. The following three, somewhat surprising, facts about \eqref{eq:rp} are also proved in \cite{FrLi}. The inequality does not hold in general for $0<\lambda<N-2$. It is strict for $N-2<\lambda<N$ unless $f\equiv 0$, but for $\lambda=N-2$ the left side can also vanish for non-trivial $f$.

A simple corollary is that  if $g$ is any other function with support in $H$ then
\begin{align}\label{eq:rpcor}
 I_\lambda[\theta g,f]
\leq \sqrt{I_\lambda[\theta f,f] \ I_\lambda[\theta g,g]} 
 \leq \frac12 I_\lambda[\theta f,f] +\frac12 I_\lambda[\theta g,g] \, .
\end{align}

The physical content of \eqref{eq:rp} for $N=3$ and $\lambda=1$ is that the interaction of an electric charge distribution with its (opposite charge) mirror image in a reflecting (Dirichlet) plane is always negative. It turns out also to be monotonic with respect to the separation distance, which implies that the charge is always drawn to the plane. (The monotonicity -- and even log-convexity -- can be deduced from \eqref{eq:rpcor} by considering various reflection planes.)

At this point, it is useful to introduce the group $\mathcal C$ of conformal transformations of $\R^N\cup\{\infty\}$ into $\R^N\cup\{\infty\}$ generated by the Euclidean group (translations and rotations) together with scaling and inversion in the unit sphere centered at the origin. Reflections through planes are in $\mathcal C$. The kernel $|x-y|^{-\lambda}$ appearing in \eqref{eq:rp} is invariant under the action of this group, except for a factor of the form $\alpha(x)\alpha(y)$. The transformation of the integral in \eqref{eq:rp} also introduces a Jacobian of the form $\beta(x)\beta(y)$. The product $\alpha(x)\beta(x)$ can be absorbed into the function $f(x)$. To be more precise, if $\gamma$ is an element of $\mathcal C$ and $f$ is a given function we define the $\lambda$-dependent transformation
\begin{equation}
 \label{eq:conftrafo}
F(x) := \left|\mathcal J_{\gamma^{-1}}(x)\right|^{(2N-\lambda)/2N} f(\gamma^{-1} x) \,.
\end{equation}
Then
\begin{align}
\label{eq:confinvintro}
I_\lambda[F,F] = \iint_{\R^N\times\R^N} \!\!\frac{\overline{F(y)}\ F(x)} {|y-x|^{\lambda}} \,dx\,dy
= \iint_{\R^N\times\R^N} \!\!\frac{\overline{f(y')}\ f(x')}{|y'-x'|^{\lambda}} \,dx'\,dy' = I_\lambda[f,f] \,.
\end{align}
Because of the group property of $\mathcal C$ it is only necessary to check this formula for translation, rotation, scaling and inversion; see \cite[Secs. 4.4 and 4.5]{LiLo}.

By a conformal transformation one can also transform the half-space into a ball. The reflection through planes becomes \textbf{reflection through spheres}. By \eqref{eq:rp} this reflection is reflection positive. This time, however, we have to supplement the geometric inversion in the sphere, $\gamma$, by a Jacobian factor $\left|\mathcal J_{\gamma^{-1}}(x)\right|^{(2N-\lambda)/2N}$. More precisely, given a ball $B$ in $\R^N$, which we may take to be centered at zero and of radius $r$, and a function $f$ with support in $B$, we define a function $\theta f$ with support in $\theta B$, the complement of $B$, by 
\begin{equation}
 \label{eq:inversion}
\theta f(x) = \left(\frac r{|x|}\right)^{2N-\lambda} f\left(\frac{r^2 x}{|x|^2}\right) \,.
\end{equation}
With this definition of $\theta f$, which now depends on $\lambda$, $N$ as well as the radius $r$, \emph{inequality \eqref{eq:rp} continues to hold}. This route to inversion positivity in spheres was pointed out to us by E. Carlen, to whom we are most grateful. Originally, we had a hands-on proof using Gegenbauer polynomials which we report here. We thank R. Askey for giving us references to \cite{Ge,As} where the necessary positivity statements are proved. 

A natural question is whether this kind of reflection positivity through planes in $\R^N$ can be generalized to \textbf{reflections through equators} in $\Sph^{N}$. The kernel is still $|s-t|^{-\lambda}$, where $s,t$ are unit vectors in $\R^{N+1}$. If we think of the sphere as embedded in $\R^{N+1}$ and use the result in \eqref{eq:rp} in $\R^{N+1}$, then the answer is immediately seen to be positive provided $N-1\leq \lambda<N$. But this is not the right way to look at it! A better way is to note that the kernel $|s-t|^{-\lambda}$ has another conformal covariance, namely under the stereographic projection $\mathcal S$ from $\Sph^N$ to $\R^N$. Under stereographic projection, reflection through equators corresponds to reflection through spheres. The dimension of the manifold is preserved and we get the correct condition $N-2\leq \lambda<N$ by using the previously obtained result.

\subsection{The HLS inequality}

The Hardy--Littlewood--Sobolev inequality for functions $f$ and $g$ on $\R^N$,
\begin{equation} \label{eq:hls1}
\Big| \, I_\lambda[f,g] \, \Big|
 \leq \mathcal{H}_{N,\lambda, p,q} \|f\|_p\|g\|_q \,,
\end{equation}
holds for all $0<\lambda <N$ and $p,q>1$ with $1/p+1/q +\lambda/N =2$ \cite[Thm. 4.3]{LiLo}. The sharp value of $\mathcal{H}_{N,\lambda, p,q}$ is known only in the \emph{diagonal case} $p=q = 2N/(2N-\lambda)$ \cite{Li}. The optimizers of \eqref{eq:hls1} are precisely the functions $f(x)=c(b^2+|x-a|^2)^{-(2N-\lambda)/2}$, $g(x)=c'f(x)$ mentioned above, where $a\in\R^N$, $b>0$ and $0\neq c,c'\in\C$. Our aim here is to prove this fact in the diagonal case when $N-2\leq \lambda<N$ by using reflection positivity through spheres instead of symmetric decreasing rearrangement, as in the original proof \cite{Li} and in \cite{CaLo}. (Recently, Carlen, Carrillo and Loss \cite{CaCaLo} have found a proof of the sharp inequality \eqref{eq:hls1} for $\lambda=N-2$ that does not use rearrangements.) Our attack on the problem will reveal the geometric significance of this class of functions, as discovered by Li and Zhu \cite{LZh}. Symmetric decreasing rearrangement is a non-linear operation whereas our reflection positivity argument is essentially linear. 

Among the diagonal cases, an important example is $\lambda= N-2$, where the kernel is 
Newton's gravitation potential. Mathematically, this case is dual to the ordinary
Sobolev inequality for $N\geq 3$, \cite[Thm. 8.3]{LiLo} $\Vert \nabla f\Vert^2_2 \geq S_N \Vert f \Vert^2_{2N/(N-2)}$, and thus the sharp constant for one gives a sharp  constant for the other.
Completely different proofs have been given for this special case \cite{Ro,Au,Ta,CoNaVi,BoLe}. Similarly, $\lambda=N-2s$ corresponds to the Sobolev inequality for $(-\Delta)^s$ when $N> 2s$.

We shall also be interested in the limiting case $\lambda\to 0$ of \eqref{eq:hls1}. Note that in the diagonal case $p=q=2N/(2N-\lambda)$ tends to $1$ in this limit, so that for non-negative functions $f$ and $g$ inequality \eqref{eq:hls1} becomes an equality. Differentiating at the end point one arrives at the logarithmic Hardy--Littlewood--Sobolev inequality
\begin{equation}\label{eq:loghls}
J[f,g] \leq \mathcal H_N
\end{equation}
for non-negative $f,g$ with $\int f \,dx=\int g \,dx=1$, where
\begin{align*}
J[f,g] := \iint_{\R^N\times\R^N} \overline{f(x)} \log\frac{1}{|x-y|} g(y) \,dx\,dy 
& - \frac{1}{2N} \int_{\R^N} f(x)\log f(x) \,dx \\ 
& - \frac{1}{2N} \int_{\R^N} g(x)\log g(x) \,dx \,.
\end{align*}
In this way one obtains the sharp constant
$$
\mathcal H_N = \frac{d}{d\lambda} \mathcal{H}_{N,\lambda, \tfrac{2N}{2N-\lambda},\tfrac{2N}{2N-\lambda}} \Big|_{\lambda=0}
$$
in \eqref{eq:loghls} from the sharp constant in \eqref{eq:hls1}. The characterization of optimizers of \eqref{eq:loghls}, however, is lost in this limit and requires additional arguments. It was shown by Carlen and Loss \cite{CaLo1} and by Beckner \cite{Be} that the optimizers are precisely the functions $f(x)=g(x)=c(b^2+|x-a|^2)^{-N}$ mentioned above, where $a\in\R^N$, $b>0$ and $c>0$ are such that the integral is equal to one. In this paper we will give a new proof of this fact for $N=1$ and $N=2$ by using reflection positivity.

In a similar way in which \eqref{eq:hls1} for $\lambda=N-2$ is equivalent to the Sobolev inequality, the logarithmic Hardy--Littlewood--Sobolev inequality \eqref{eq:loghls} for $N=2$ is equivalent to Onofri's inequality and for $N=1$ to the Lebedev--Milin inequality, see \cite{CaLo1,Be}. For alternative proofs of Onofri's inequality and its generalizations we refer to \cite{On,Ho,OsPhSa,Wi,ChYa}.


\section{Main results}

We shall prove

\begin{theorem}[\textbf{HLS inequality}]\label{mainhls}
 Let $0<\lambda<N$ if $N=1,2$ and $N-2\leq\lambda<N$ if $N\geq 3$. If $p=q=2N/(2N-\lambda)$, then \eqref{eq:hls1} holds with
\begin{equation}
 \label{eq:hlsconst}
\mathcal H_{N,\lambda,p,p} = \pi^{\lambda/2} \frac{\Gamma((N-\lambda)/2)}{\Gamma(N-\lambda/2)} \left(\frac{\Gamma(N)}{\Gamma(N/2)}\right)^{1-\lambda/N} \,.
\end{equation}
Equality holds if and only if $f\equiv 0$ or $g\equiv 0$ or
$$
f(x) = c \left(b^2 +|x-a|^2\right)^{-(2N-\lambda)/2}
\quad\text{and}\quad
g(x) = c' \left(b^2 +|x-a|^2\right)^{-(2N-\lambda)/2} \,,
$$
for some $a\in\R^N$, $b>0$ and $c,c'\in\C$.
\end{theorem}

Our second main result is

\begin{theorem}[\textbf{Logarithmic HLS inequality}]\label{mainloghls}
 If $N=1,2$, then \eqref{eq:loghls} holds with
\begin{equation}
 \label{eq:loghlsconst}
\mathcal H_{N} = 
\frac12 \log \pi + \frac12(\psi(N)-\psi(N/2)) - \frac1N \log\frac{\Gamma(N)}{\Gamma(N/2)}
\end{equation}
for any non-negative functions $f,g$ on $\R^N$ satisfying
$$
\int_{\R^N} f(x)\,dx= \int_{\R^N} g(x) \,dx =1
$$
and $\int f(x) \log_+ ( f(x) (\frac{1+|x|^2}{2})^N ) \,dx<\infty$, $\int g(x) \log_+ ( g(x) (\frac{1+|x|^2}{2})^N ) \,dx<\infty$. Here $\psi=(\log\Gamma)'$ is the digamma function. Equality holds if and only if
$$
f(x) =g(x) = \frac{2^{N-1} \ \Gamma((N+1)/2)}{\pi^{(N+1)/2}} \ b^{N} \left(b^2 +|x-a|^2\right)^{-N}
$$
for some $a\in\R^N$ and $b>0$.
\end{theorem}

Since $\psi(1)=-\gamma$, $\psi(2)=1-\gamma$ and $\psi(1/2)=-2\log 2-\gamma$ ($\gamma$ the Euler-Mascheroni constant), one finds
$$
\mathcal H_1= \log (2\pi)\,,
\qquad
\mathcal H_2= \frac12 (1+\log \pi) \,.
$$


\subsection*{Outline of the proofs of Theorems \ref{mainhls} and \ref{mainloghls}}

As observed in \cite{Li,CaLo,CaLo1}, the functionals $I_\lambda$ and $J$ are conformally invariant. This implies in particular that the values of $I_\lambda[f]:=I_\lambda[f,f]$ and $J[f]:=J[f,f]$ do not change if $f$ is inverted on the surface of a ball or reflected on a hyperplane. To state this property precisely, we need to introduce some notation. The dependence on the fixed parameter $0\leq \lambda<N$ will not be reflected in the notation.

Let $B=\{x \in\R^N : \ |x-a|<r\}$, $a\in\R^N$, $r>0$, be an open ball and denote by
$$
\Theta_B(x) := \frac{r^2(x-a)}{|x-a|^2} +a
$$
the inversion of a point $x\neq a$ through the boundary of $B$. This map on $\R^N$ can be lifted to an operator acting on functions $f$ on $\R^N$ according to
$$
(\Theta_B f)(x) := \left(\frac{r}{|x-a|}\right)^{2N-\lambda} f(\Theta_B(x)) \,.
$$
(Strictly speaking, $\Theta_B f$ is not defined at the point $x=a$.) Note that both the map and the operator $\Theta_B$ satisfy $\Theta_B^2=I$, the identity. The crucial property for us is that
\begin{equation}
 \label{eq:confinv}
I_\lambda[f]=I_\lambda[\Theta_B f]
\end{equation}
if $\lambda>0$ and $J[f]= J[\Theta_B f]$ for $\lambda=0$.

Similarly, let $H=\{x \in\R^N : \ x\cdot e>t \}$, $e\in\Sph^{N-1}$, $t\in\R$, be a half-space and denote by
$$
\Theta_H(x) := x+ 2(t-x\cdot e)
$$
the reflection of a point $x$ on the boundary of $H$. The corresponding operator is defined by
$$
(\Theta_H f)(x) := f(\Theta_H(x))
$$
and it again satisfies $\Theta_H^2=I$. Moreover,
\begin{equation}
 \label{eq:confinvh}
I_\lambda[f]=I_\lambda[\Theta_H f] \,.
\end{equation}
if $\lambda>0$ and $J[f]= J[\Theta_H f]$ for $\lambda=0$. Our first ingredient in the proof of Theorem~\ref{mainhls} is the following.

\begin{theorem}[\textbf{Reflection positivity in planes and spheres}]\label{infrefpos}
 Let  $0\leq\lambda<1$ if $N=1$, $N-2\leq\lambda<N$ if $N\geq 2$ and let $B\subset\R^N$ be either a ball or a half-space. For $f\in L^{2N/(2N-\lambda)}(\R^N)$ define
\begin{equation*}
 f^i (x) :=
\begin{cases}
 f(x) & \text{if}\ x\in B\,,\\
\Theta_B f(x) & \text{if}\ x\in \R^N\setminus B\,,
\end{cases}
\qquad
f^o (x) :=
\begin{cases}
 \Theta_B f(x) & \text{if}\ x\in B \,, \\
f(x) & \text{if}\ x\in \R^N\setminus B \,.
\end{cases}
\end{equation*}
Then for $\lambda>0$ one has
\begin{equation}\label{eq:infrefpos}
\frac12 \left(I_\lambda[f^i]+I_\lambda[f^o]\right) \geq I_\lambda[f] \,,
\end{equation}
and for $\lambda=0$ under the additional assumptions that $f\geq 0$ and $\int_B f \,dx = \int_{B^c} f\,dx$ one has
\begin{equation}\label{eq:infrefposlog}
 \frac12 \left(J[f^i]+J[f^o]\right) \geq J[f] \,.
\end{equation}
If $\lambda>N-2$, then inequalities \eqref{eq:infrefpos} and \eqref{eq:infrefposlog} are strict unless $f=\Theta_B f$.
\end{theorem}

For half-spaces and $\lambda=N-2$ (the Newtonian case) this theorem was long known to quantum field theorists \cite{GlJa,Li0,OsSc}. The half-space case with $N-2<\lambda<N$ (but not the strictness for $\lambda>N-2$) was apparently first proved by Lopes and Mari\c{s} \cite{LoMa}. The case of balls seems to be new for all $\lambda$.

Our second main ingredient is the following generalization of a theorem of Li and Zhu \cite{LZh}; see also \cite{L}. We refer to \cite{FrLi} for the proof.

\begin{theorem}[\textbf{Characterization of inversion invariant measures}]\label{yyl}
 Let $\mu$ be a finite, non-negative measure on $\R^N$. Assume that 
\begin{itemize}
 \item[(A)] \label{ass:yyl}
for any $a\in\R^N$ there is an open ball $B$ centered at $a$ and for any $e\in\Sph^{N-1}$ there is an open half-space $H$ with interior unit normal $e$ such that
\begin{equation}
 \label{eq:yyl}
\mu(\Theta^{-1}_B(A))=\mu(\Theta^{-1}_H(A))=\mu(A)
\qquad \text{for any Borel set}\ A\subset\R^N \,.
\end{equation}
\end{itemize}
Then $\mu$ is absolutely continuous with respect to Lebesgue measure and
$$
d\mu(x) = c \left(b^2 +|x-a|^2\right)^{-N} dx
$$
for some $a\in\R^N$, $b>0$ and $c\geq 0$.
\end{theorem}

We emphasize that $B$ and $H$ in assumption (A) divide $\mu$ in half, in the sense that $\mu(B)=\mu(\R^N\setminus\overline B)$ and $\mu(H)=\mu(\R^N\setminus\overline H)$. By a change of variables one finds that for absolutely continuous measures $d\mu=v \,dx$ assumption (A) is equivalent to the fact that for any $a\in\R^N$ there is an $r_a>0$ and a set of full measure in $\R^N$ such that for any $x$ in this set
\begin{equation}\label{eq:pweq}
v(x) = \left(\frac{r_a}{|x-a|}\right)^{2N} v\left( \frac{r_a^2 (x-a)}{|x-a|^2}+a \right) \,,
\end{equation}
and similarly for reflections. We emphasize that the assumption that $\mu$ is finite is essential in Theorem \ref{yyl}, since $d\mu(x)=|x|^{-2N}dx$ also satisfies assumption (A).

Theorems \ref{mainhls} and \ref{mainloghls} follow from Theorems \ref{infrefpos} and \ref{yyl}. Since we have shown this in \cite{FrLi} for the Hardy--Littlewood--Sobolev inequality, we concentrate here on its logarithmic version.

\begin{proof}
 We first note that we can restrict our attention to the case $f=g$ because $J[\frac12(f+g)] \geq J[f,g]$ with strict inequality unless $f=g$. To verify this claim, we put
\begin{equation}\label{eq:i0}
I_0[f,g] := \iint_{\R^N\times\R^N} \overline{f(x)}\ \log\frac{1}{|x-y|} \ g(y) \,dx\,dy
\end{equation}
and $I_0[f]:=I_0[f,f]$. Now for real $f$ and $g$, $I_0[\frac12(f+g)]= I_0[f,g]+I_0[\tfrac12(f-g)]$ and we shall see in the proof of Lemma \ref{repr} that $I_0[h]\geq 0$ for all $h$ with $\int h\,dx=0$. Hence $I_0[\frac12(f+g)]\geq I_0[f,g]$. Moreover, since $x\log x$ is a strictly convex function of $x>0$,
$$
 \int \frac12(f+g) \log \frac12(f+g) \,dx \leq  \frac12 \left( \int f\log f\,dx + \int g\log g\,dx \right) 
$$
with strict inequality unless $f\equiv g$. This proves the claim.

Next, we claim that the supremum
$$
\sup\left\{ J[h] :\ h\geq 0 \,, \int h\,dx=1 \, \int h(x) \log_+ \left( h(x) \left(\tfrac12\left(1+|x|^2\right)\right)^N \right) \,dx<\infty  \right\}
$$
is attained and given by the right side of \eqref{eq:loghlsconst}. Indeed, differentiating the Hardy--Littlewood--Sobolev inequality \eqref{eq:hls1} at the endpoint $\lambda=0$ we see that the right side of \eqref{eq:loghlsconst} is an upper bound for $J[h]$. On the other hand, for the $f$ given in Theorem \ref{mainloghls} one can compute that $J[f]$ is given by the right side of \eqref{eq:loghlsconst}. This proves the statement about the sharp constant and we are left with characterizing the optimizers.

 Let $f$ be an optimizer, that is, a non-negative function $f$ with $\int f\,dx=1$ for which the above supremum is attained. For any point $a$ there is a ball $B$ centered at $a$ such that $\int_{B} f \,dx = \int_{\R^N\setminus B} f \,dx$. We note that if $f^i$ and $f^o$ are defined as in Theorem \ref{infrefpos}, then $\|f^i\|_1=\|f^o\|_1=\|f\|_1=1$. Moreover, by \eqref{eq:infrefposlog}, $\tfrac12( J[f^i]+J[f^o] )\geq J[f]$ and hence, in particular, $\max\{J[f^i], J[f^o]\} \geq J[f]$. By the maximizing property of $f$ this inequality cannot be strict, and therefore we conclude that $J[f^i]= J[f^o]= J[f]$, that is, both $f^i$ and $f^o$ are optimizers as well.

 In order to continue the argument we assume first that $N=1$. Since we have just shown that one has equality in \eqref{eq:infrefposlog}, the second part of Theorem \ref{infrefpos} implies that $f=\Theta_{B} f$. By a similar argument one deduces that $f=\Theta_H f$ for any half-space such that $\int_H f \,dx = \int_{\R^N\setminus H} f \,dx$. Therefore the measure $f\,dx$ satisfies the assumption of Theorem \ref{yyl}, and hence $f$ has the form claimed in Theorem \ref{mainloghls}.

Now assume that $N=2$. The difference from the previous case is that there is no strictness assertion in Theorem \ref{infrefpos} (indeed, equality in \eqref{eq:infrefpos} can hold without having $f=\Theta_B f$), so we need an additional argument in the spirit of \cite{Lo} to conclude that $f=\Theta_B f$ for any ball and half-space with $\int_B f \,dx = \int_{\R^N \setminus B} f \,dx$.

Before proceeding we shall show that $f$ (and therefore also $f^o$ and $f^i$) are a.e. positive. Indeed, if $f$ would vanish on a set $K$ of positive (but finite) measure we could take $f_\epsilon := (f+\epsilon\chi_K)/(1+\epsilon |K|)$ as a trial function and find
$$
J[f_\epsilon]= J[f] - \frac12 |K| \epsilon\log\epsilon + \mathcal O(\epsilon)
$$
as $\epsilon\to 0$. This contradicts the maximizing property of $f$ and shows that $|K|=0$.

In the first part of the proof we have seen that $f^o$ (and $f^i$) are optimizers. Using that they are positive a.e. we find that they satisfy the Euler-Lagrange equations
$$
\int_{\R^2} f(y) \log\frac{1}{|x-y|} \,dy - \frac14 \log f(x) = \mu \,,
\qquad
\int_{\R^2} f^o(y) \log\frac{1}{|x-y|} \,dy - \frac14 \log f^o(x) = \mu \,.
$$
Here the Lagrange multipliers coincide since $J[f]=J[f^o]$ and $\|f^o\|_1=\|f\|_1=1$. The functions $u:=\log (8\pi f)$ and $u^o := \log (8\pi f^o)$ satisfy the equations
$$
 -\Delta u = e^u \,,
\qquad
-\Delta u^o = e^{u^o} \,.
$$
Since $e^u, e^{u^o}\in L^1(\R^2)$, we deduce from \cite{BrMe} that $u, u^o\in L^\infty(\R^2)$. The function $w:=u-u^o$ satisfies $-\Delta w+Vw=0$ with
$$
V(x):= - \frac{e^{u(x)} - e^{u^o(x)}}{u(x)-u^o(x)} 
= - \int_0^1 e^{t u(x)+ (1-t)u^o(x)} \,dt \,.
$$
Since $u$ and $u^o$ are bounded, $V$ is so as well. Since $w\equiv 0$ in $\R^N\setminus B$, the unique continuation theorem \cite{Ca} implies that $w\equiv 0$ everywhere. Hence $f=\Theta_Bf$ and we can deduce Theorem \ref{mainloghls} again from Theorem \ref{yyl}.
\end{proof}


\section{Reflection positivity in planes and spheres}\label{sec:pos}

Our goal in this section is to prove Theorem \ref{infrefpos}. In Subsections \ref{sec:repr} and \ref{sec:reprball} we consider the cases of half-spaces and balls, respectively, and derive representation formulas for $I_\lambda[\Theta_H f, f]$ and $I_\lambda[\Theta_B f, f]$. In Subsection \ref{sec:infrefposproof} we prove Theorem \ref{infrefpos}.


\subsection{Reflection positivity in planes}
\label{sec:repr}

Throughout this subsection we assume that $H=\{x\in\R^N:\ x_N> 0 \}$. The key for proving Theorem \ref{infrefpos} is the following explicit formula for $I_\lambda[\Theta_H f,f]$. Recall that $I_0[f,g]$ was defined in \eqref{eq:i0}.

\begin{lemma}[\textbf{Representation formula}]
 \label{repr}
 Let $0\leq \lambda<1$ if $N=1$ and $N-2\leq \lambda<N$ if $N\geq 2$. Let $f\in L^{2N/(2N-\lambda)}(\R^N)$ be a function with support in $\overline H = \{x\in\R^N:\ x_N\geq 0 \}$. If $\lambda=0$ assume, in addition, that $\int_{\R^N} f(x)\,dx=0$.
\begin{enumerate}
 \item If $\lambda>N-2$, then
\begin{equation}
 \label{eq:repr}
I_\lambda[\Theta_H f,f] = c_{N,\lambda} 
\int_{\R^{N-1}} d\xi'
\int_{|\xi'|}^\infty d\tau \frac{\tau^2}{(\tau^2-|\xi'|^2)^{(N-\lambda)/2}}
\left| \int_\R \frac{\hat f(\xi)} {\tau^2+\xi_N^2} \,d\xi_N \right|^2
\end{equation}
where
\begin{equation*}
 c_{N,\lambda} =
\begin{cases}
 2^{N+1-\lambda} \pi^{(N-4)/2} \, \frac{\sin(\pi(N-\lambda)/2) \ \Gamma((N-\lambda)/2)}{\Gamma(\lambda/2)} 
& \text{if}\ \lambda>0 \,, \\
\frac2\pi & \text{if}\ \lambda=0 \,.
\end{cases}
\end{equation*}
\item 
If $\lambda=N-2$, then
\begin{equation}
 \label{eq:reprcou}
I_{N-2}[\Theta_H f,f] = c_{N,N-2} \ 
\int_{\R^{N-1}} d\xi' |\xi'|
\left| \int_\R \frac{\hat f(\xi)} {|\xi'|^2+\xi_N^2} \,d\xi_N \right|^2 \,.
\end{equation}
where
\begin{equation*}
 c_{N,N-2} =
\begin{cases}
2 & \text{if}\ N=2 \,, \\
\frac{4 \pi^{(N-2)/2}}{\Gamma((N-2)/2)} & \text{if}\ N\geq 3 \,.
\end{cases}
\end{equation*}
\end{enumerate}
In any case, $c_{N,\lambda}>0$.
\end{lemma}

When $N=1$, we use the convention that $\R^{N-1}=\{0\}$ and that $d\xi'$ gives measure $1$ to this point.

The crucial point of Lemma \ref{repr} is, of course, that the right sides of \eqref{eq:repr} and \eqref{eq:reprcou} are non-negative. Indeed, in Subsection \ref{sec:infrefposproof} we shall see that the right side of \eqref{eq:repr} is strictly positive unless $f\equiv 0$. We refer to \cite{FrLi} for the facts that $I_\lambda[\Theta_H f,f]$ is not necessarily non-negative for $\lambda<N-2$ and that $I_{N-2}[\Theta_H f,f]$ can vanish even if $f\not\equiv 0$.

\begin{proof}
 For $N=1$ and $0<\lambda<1$ one has
$$
I_\lambda[\Theta_H f,f] = \int_0^\infty \int_0^\infty \frac{\overline{f(x)}\ f(y)}{(x+y)^{\lambda}} \,dx\,dy 
= \frac1{\Gamma(\lambda)} \int_0^\infty \frac{d\tau}{\tau^{1-\lambda}} \left| \int_0^\infty e^{-\tau x} f(x) \,dx \right|^2 \,.
$$
Recalling that $f(x)=0$ for $x\leq 0$ and using that $e^{-\tau|x|}$ has Fourier transform $(2/\pi)^{1/2} \tau/(\xi^2+\tau^2)$ we can write
\begin{equation}
 \label{eq:cauchy}
\int_0^\infty e^{-\tau x} f(x) \,dx = \sqrt\frac2\pi \ \tau \int_\R \frac{\hat f(\xi)}{\xi^2+\tau^2} \,d\xi \,.
\end{equation}
Since $c_{1,\lambda}= 2/(\pi \Gamma(\lambda))$ we have shown the assertion in this case. If $f$ has integral zero, then according to the above
\begin{align*}
\lambda^{-1} \int_0^\infty \int_0^\infty \overline{f(x)}\ \left(\frac1{(x+y)^{\lambda}} - 1 \right) f(y) \,dx\,dy
& = \lambda^{-1} I_\lambda[f,\Theta_H f] \\ 
& = \frac1{\lambda\Gamma(\lambda)} \int_0^\infty \frac{d\tau}{\tau^{1-\lambda}} \left| \int_0^\infty e^{-\tau x} f(x) \,dx \right|^2 \,.
\end{align*}
Letting $\lambda\to 0$ from above and noting that $\lambda\Gamma(\lambda)\to 1$ we obtain
$$
I_0[\Theta_H f,f] = \int_0^\infty \frac{d\tau}{\tau} \left| \int_0^\infty e^{-\tau x} f(x) \,dx \right|^2 \,.
$$
This together with \eqref{eq:cauchy} proves the claim for $N=1$ and $\lambda=0$.

Finally, we prove the assertion for $N=2$ and $\lambda=0$, the other cases being contained in \cite{FrLi}. Since $f$ has integral zero, we have
\begin{align*}
I_0[\Theta_H f,f] 
& = \int_H \int_H \overline{f(x)} \log\frac1{\sqrt{(x_1-y_1)^2+(x_2+y_2)^2}} f(y) \,dx\,dy \\ 
& = (2\pi)^{-1} \int_{\R^2} 
\int_H \int_H \overline{f(x)} \ \frac{e^{i\xi_1(x_1-y_1) + i\xi_2(x_2+y_2)}}{|\xi|^2} f(y) \,dx\,dy\,d\xi \\
& = \int_{\R} K_{\xi_1}[F_{\xi_1}] \,d\xi_1 \,,
\end{align*}
where
$$
F_{\xi_1}(t) := (2\pi)^{-1/2} \int_{\R} f(x_1,t) e^{-i\xi_1x_1} \,dx_1 \,,
$$
and
$$
K_{\xi_1}[\phi] = \int_0^\infty \int_0^\infty \overline{\phi(t)} k_{\xi_1}(t+s) \phi(s) \,ds\,dt \,,
\quad
k_{\xi_1}(t) :=  \int_\R \frac{e^{i\xi_2 t} }{\xi_1^2+\xi_2^2} \,d\xi_2 \,.
$$
By the residue theorem $k_{\xi_1}(t) = \pi |\xi_1|^{-1} e^{-t|\xi_1|}$, and hence
$$
K_{\xi_1}[\phi] = \pi |\xi_1|^{-1} \left| \int_0^\infty \! e^{-t|\xi_1|} \phi(t) \,dt \right|^2 \,.
$$
In view of \eqref{eq:cauchy} this is the claimed formula.
\end{proof}


\subsection{Reflection positivity in spheres}\label{sec:reprball}

Throughout this subsection we assume that $B=\{x\in\R^N:\ |x|<1 \}$ is the unit ball in $\R^N$. For $N\geq 2$ we denote by $\{Y_{l,m}\}$ an orthonormal basis of $L^2(\Sph^{N-1})$, where $Y_{l,m}$ is a spherical harmonic of degree $l$. The index $l$ runs through $\N_0$ and, for any fixed $l$, $m$ runs through a certain set of $l$-dependent cardinality. For $N=1$ we have $\Sph^0=\{-1,1\}$ and we put $Y_{l,0}(\omega) := 2^{-1/2} \omega^l$ for $l=0,1$. Here $l$ assumes only the values $0$ and $1$, and $m=0$. Furthermore, we shall need the Gegenbauer polynomials $C_k^{(\alpha)}$. 
Their definition as well as basic properties may be found, e.g., in \cite{GrRy}. The constants
\begin{equation}\label{eq:gegconst}
c^{\alpha,\beta}_{k,l} := \int_{-1}^1 C_k^{(\alpha)}(\tau) \ C_l^{(\beta)}(\tau) \ (1-\tau^2)^{\beta-1/2} \,d\tau
\end{equation}
for $\alpha,\beta>- 1/2$ will appear in our analysis below.

\begin{lemma}[\textbf{Representation formula}]
 \label{reprball}
 Let $0\leq \lambda<N$ and let $f\in L^{2N/(2N-\lambda)}(\R^N)$ be a function with support in $\overline B = \{x\in\R^N:\ |x|\leq 1 \}$. If $\lambda=0$ assume, in addition, that $\int_{\R^N} f(x)\,dx=0$. Then
\begin{equation}
 \label{eq:reprball}
I_\lambda[\Theta_B f, f] 
= \sum_{l,m} \sum_{k=0}^\infty \kappa_{N,l}  \ d_{k,l}^{(\lambda)}
\left| \int_B \overline{Y_{l,m}(x/|x|)} f(x)  |x|^{k} \,dx \right|^2 \,.
\end{equation}
If $N\geq 2$, then
\begin{equation*}
d_{k,l}^{(\lambda)} := 
\begin{cases}
c_{k,l}^{\lambda/2,(N-2)/2} & \text{if} \ \lambda>0\,,\\
\frac12  c_{k,l}^{0,(N-2)/2} & \text{if} \ \lambda=0 \,,
\end{cases}
\end{equation*}
and if $N=1$, then $d_{0,l}^{(\lambda)} := \left( 1+(-1)^{k+l} \right)$ for $k=0$ and
\begin{equation*}
d_{k,l}^{(\lambda)} := 
\begin{cases}
\frac{\lambda \cdots (\lambda+k-1)}{k!} \left( 1+(-1)^{k+l} \right) & \text{if} \ \lambda>0\,,\\
\frac1k \left( 1+(-1)^{k+l} \right) & \text{if} \ \lambda=0 \,,
\end{cases}
\end{equation*}
for $k\geq 1$. Moreover,
\begin{equation}
 \label{eq:kappa}
\kappa_{N,l} := 
\begin{cases}
1 & \text{if} \ N=1\,,\, l=0,1\,,\\
 2 & \text{if} \ N=2\,,\, l=0 \,,\\
l & \text{if} \ N=2\,,\, l\geq 1 \,,\\
(4\pi)^{\frac{N-2}2} \ \dfrac{l! \ \Gamma((N-2)/2)}{(l+N-3)!} & \text{if} \ N\geq 3\,,\, l\geq 0 \,.
\end{cases}
\end{equation}
\end{lemma}

\begin{proof}
We begin with the case $\lambda>0$ and note that
 $$
I_\lambda[\Theta_B f, f]
= \int_{B^c} \,dx \int_{B}\,dy \frac{\overline{f(x/|x|^2)}\ f(y)}{|x|^{2N-\lambda} |x-y|^\lambda }
= \iint_{B\times B} \frac{\overline{f(x)}\ f(y)}{\left|\left|x\right|^{-1} x - \left|x\right| y \right|^\lambda} \,dx\,dy \,.
$$
Next, we write $x=r\omega$ with $\omega\in\Sph^{N-1}$ and decompose $f$ into spherical harmonics
$$
f(x) = \sum_{l,m} f_{l,m}(r) Y_{l,m}(\omega) \,,
\qquad f_{l,m}(r) := \int_{\Sph^{N-1}} \overline{Y_{l,m}(\omega)} f(r\omega) \,d\omega \,.
$$
We shall see that
\begin{equation}
 \label{eq:adecomp}
\iint_{B\times B} \frac{\overline{f(x)}\ f(y)}{\left|\left|x\right|^{-1} x - \left|x\right| y \right|^\lambda} \,dx\,dy
= \sum_{l,m} \kappa_{N,l} \, A_l[f_{l,m}]
\end{equation}
where the functionals $A_l$ are of the form
$$
A_l[\phi] = \int_0^1 \int_0^1 \overline{\phi(r)} \ a_l(rs) \ \phi(s) \, r^{N-1}dr \, s^{N-1}ds
$$
with functions $a_l$ to be defined below.

\emph{Case $N=1$}. An easy calculation shows that
$$
a_l(r) = \frac{1}{(1-r)^\lambda} + (-1)^l \frac{1}{(1+r)^\lambda}
$$
for $l=0$ and $l=1$. Writing
$$
(1-r)^{-\lambda} = \sum_{k=0}^\infty \gamma_{\lambda,k} r^k
\qquad \text{with}\ \gamma_{\lambda,0}=1 \,, \ \
\gamma_{\lambda,k}=\frac{\lambda \cdots \left(\lambda+k-1\right)}{k!} \ \text{for}\ k\geq 1 \,,
$$
we deduce that
$$
A_l[\phi] = \sum_{k=0}^\infty \left( 1 + (-1)^{l+k}  \right) \gamma_{\lambda,k} \left| \int_0^1 \phi(r) r^k \,dr \right|^2 \,.
$$

\emph{Case $N\geq 2$}. We recall the Funk-Hecke formula (see, e.g., \cite[Sec. 11.4]{ErMaObTr}), which states that for $K\in L^1((-1,1),(1-t^2)^{(N-3)/2}dt)$ and any spherical harmonic $Y$ of degree $l$ on $\Sph^{N-1}$ one has
\begin{equation*}
 \int_{\Sph^{N-1}} K(\omega\cdot\omega') Y(\omega') \,d\omega' = 
\kappa_{N,l} \left( \int_{-1}^1 K(\tau) C_l^{(\frac{N-2}2)}(\tau) (1-\tau^2)^{\frac{N-3}2} \,d\tau \right) Y(\omega)
\end{equation*}
for all $\omega\in\Sph^{N-1}$ with $\kappa_{N,l}$ given in \eqref{eq:kappa}. This implies that \eqref{eq:adecomp} holds with
$$
a_l(r) = \int_{-1}^1 (1-2r\tau +r^2 )^{-\lambda/2} C_l^{(\frac{N-2}2)}(\tau) (1-\tau^2)^{\frac{N-3}2} \,d\tau \,.
$$
Using the generating function identities for the Gegenbauer polynomials,
\begin{equation}
 \label{eq:gengegen}
\left(1 -2 r\tau  + r^{2} \right)^{-\alpha} = \sum_{k=0}^\infty C_k^{(\alpha)}(\tau) \ r^k \,\quad \text{if}\ \alpha\neq 0\,,
\end{equation}
we find that
$$
a_l(r) = \sum_{k=0}^\infty c_{k,l}^{\lambda/2,(N-2)/2} \, r^k \,,
$$
and therefore
$$
A_l[\phi] = \sum_{k=0}^\infty c_{k,l}^{\lambda/2,(N-2)/2} \left| \int_0^1 \phi(r) r^{k+N-1}dr \right|^2 \,.
$$
This proves the assertion for $\lambda>0$. The proof for $\lambda=0$ is similar, replacing \eqref{eq:gengegen} by
\begin{equation*}
1-\log \left(1 -2 r\tau  + r^{2} \right) = \sum_{k=0}^\infty C_k^{(0)}(\tau) \ r^k \,.
\qedhere
\end{equation*}
\end{proof}

We emphasize that, in contrast to Lemma \ref{repr}, we have stated Lemma \ref{reprball} for $\lambda$ from the whole range $[0,N)$. The restriction $\lambda\geq N-2$ for $N\geq 3$ comes in when discussing the sign of the coefficients $d_{k,l}^{(\lambda)}$. Since obviously $d_{k,l}^{(\lambda)}\geq 0$ for $N=1$ and $0\leq \lambda<1$, we shall concentrate now on the case $N\geq 2$. Since $C_k^{(\alpha)}$ is an even function for even $k$ and an odd function for odd $k$, one has $c^{\alpha,\beta}_{k,l}=0$ if $k-l$ is odd. Moreover, since the $C_l^{(\beta)}$, $l=0,1,2,\ldots$, are the orthogonal polynomials with respect to the measure $(1-\tau^2)^{\beta-1/2} \,d\tau$, one has $c^{\alpha,\beta}_{k,l}=0$ if $k<l$. This leaves us with the case $k=l+2n$ for $n=0,1,2,\ldots$. In order to simplify the following discussion, we assume that $\alpha\neq 0$ and $\beta\neq 0$. (The formulas in the other cases are easily obtain from those below by using that $C_0^{(0)} = C_0^{(\alpha)}\equiv 1$ for any $\alpha$ and that $\lim_{\alpha\to 0} \alpha^{-1} C_k^{(\alpha)}(\tau) = C_k^{(0)}(\tau)$ for $k\geq 1$.) For $\alpha=\beta$ one has \cite[(7.313)]{GrRy}
\begin{equation}
 \label{eq:intaa}
c^{\alpha,\alpha}_{k,l} = 
\delta_{k,l}\ \frac{2^{1-2\alpha} \pi \ \Gamma(k+2\alpha)}{k!\ \ (k+\alpha) \ \Gamma(\alpha)^2} \,.
\end{equation}
By a theorem of Gegenbauer \cite{Ge} (see also \cite[Sec. 7]{As}) one has for $\alpha\neq\beta$
$$
c^{\alpha,\beta}_{k,l} = \frac{2^{1-2\beta} \pi \ \Gamma(l+n+\alpha) \ \Gamma(n+\alpha-\beta) \ \Gamma(l+2\beta)}{l!\ n! \ \Gamma(\alpha)\ \Gamma(\beta)\ \Gamma(\alpha-\beta)\ \Gamma(l+n+\beta+1)} \,.
$$
Hence by inspection,
$$
c^{\alpha,\beta}_{k,l} >0 \qquad \text{if}\ \alpha>\beta\geq 0 \ \text{and}\ k=l+2n \ \text{for}\ n\in\N_0 \,.
$$
To summarize, the coefficients $d_{k,l}^{(\lambda)}$ are non-negative if $N-2\leq\lambda<N$ and $N\geq 2$.


\subsection{Proof of Theorem \ref{infrefpos}}
\label{sec:infrefposproof}

We give the proof only in the case $\lambda=0$, the other cases being similar; see also \cite{FrLi}. We begin by considering a half-space, which after a translation and a rotation we may assume to be $H=\{x : \ x_N>0\}$. A simple calculation shows that
\begin{align*}
& \frac12 \left(J[f^i]+J[f^o]\right) - J[f] \\
& = \int_H \int_H (f(x)-f(x',-x_N)) \log\frac{1}{\sqrt{|x'-y'|^2+(x_N+y_N)^2}}  (f(y)-f(y',-y_N)) \,dx\,dy \,.
\end{align*}
Defining $g:=f-\Theta_H f$ in $H$ and $g:=0$ in $H^c$, the right side can be rewritten as $I_0[\Theta_H g,g]$. The assumption $\int_H f\,dx = \int_{H^c} f\,dx$ implies that $\int g\,dx =0$. Hence according to Lemma \ref{repr}, $I_0[\Theta_H g,g] \geq 0$.

Now assume that $N=1$ and $I_0[\Theta_H g,g]=0$. Then by \eqref{eq:repr} and \eqref{eq:cauchy} one has
\begin{equation}
 \label{eq:vanish}
\int_0^\infty e^{-\tau t} g(t) \,dt = 0
\qquad\text{for a.e.}\ \tau\geq 0 \,.
\end{equation}
By the uniqueness of the Laplace transform, $g\equiv 0$, which implies $f\equiv \Theta_H f$, as claimed.

Next, we consider a ball which we may assume to be $B=\{x\in\R^N:\ |x|<1 \}$. By a similar calculation as before,
\begin{align*}
& \frac12 \left(J[f^i]+J[f^o]\right) - J[f] \\
& \qquad = \int_{B^c} \,dx \int_{B} \,dy  \left(\Theta_B f(x) - f(x) \right) \log\frac{1}{|x-y|} \left(f(y)-\Theta_B f(y)\right) \\
& \qquad = I_0[\Theta_B g,g] \,,
\end{align*}
where now $g:=f-\Theta_B f$ in $B$ and $g:=0$ in $B^c$. Again, the assumption $\int_B f\,dx = \int_{B^c} f\,dx$ implies that $\int g\,dx =0$. As discussed at the end of the previous subsection, the coefficients $d_{k,l}^{(\lambda)}$ in \eqref{eq:reprball} are non-negative. Hence by Lemma \ref{reprball}, $I_0[\Theta_B g,g]\geq 0$.

Similarly as before, if $N=1$ and $I_0[\Theta_B g,g]=0$, then by \eqref{eq:reprball}
$$
\int_0^1 g_e(x)  x^{k} \,dx = 0
\qquad\text{for all even integers}\ k 
$$
and
$$
\int_0^1 g_o(x)  x^{k} \,dx = 0
\qquad\text{for all odd integers}\ k\,, 
$$
where $g_e$ and $g_o$ denote the even and odd parts of $g$. Changing variables $x^2=y$, we see that the functions $y^{-1/2} g_e(\sqrt y)$ and $g_o(\sqrt{y})$ are integrable on $[0,1]$ and their integral against any polynomial is zero. This implies $g_e\equiv g_o\equiv 0$ and hence $f=\Theta_B f$. This completes the proof.
\lanbox 



\bibliographystyle{amsalpha}

\end{document}